\documentclass[10pt,leqno]{article}
\usepackage{amsmath,amsfonts,amsthm}
\usepackage{fullpage}

\newcommand{\ds}{\displaystyle}
\newcommand{\la}{\lambda}

\newtheorem{thm*}{Theorem}
\newtheorem{thm}{Theorem}

\newtheorem{prop*}{Proposition}

\newtheorem{cor}{Corollary}

\date{}

\numberwithin{equation}{section}
\begin{document}
\title{A Branching Process for Virus Survival}

\bigskip
\author{J. Theodore Cox\thanks{ Supported in part by NSF Grant
  No. 0803517}\\ {\it Syracuse University} 
 \and
  Rinaldo B. Schinazi\\ {\it University of Colorado, Colorado Springs}}
\maketitle

{\bf Abstract.} Quasispecies theory predicts that there
is a critical mutation probability above which a viral
population will go extinct. Above this threshold the virus
loses the ability to replicate the best adapted genotype,
leading to a population composed of low replicating mutants
that is eventually doomed. We propose a new branching model
that shows that this is not necessarily so. That is, a population
composed of ever changing mutants may survive. 



\section{Introduction.}

Compared to other species an RNA virus has a very high
mutation rate and a great deal of genomic diversity.  Hence,
a virus population can be thought of as an ensemble of
related genotypes called quasispecies, see Eigen (1971) and
Eigen and Schuster (1977).  From the virus point of view a
high mutation rate is advantageous because it may create
rather diverse virus genomes, this may overwhelm the immune
system of the host and ensure survival of the virus
population, see Vignuzzi et al.\ (2006).  On the other hand,
a high mutation rate may result in many nonviable
individuals and hurt the quasispecies, see Sanjuan et al.\
(2004) and Elena and Moya (1999). It seems therefore that
mutation rates should be high but not too high. A simple
mathematical model makes this point. Consider a virus
population having genomes 1 and 2, where genome 1 has a
higher replication rate $a_1$ and genome 2 has a lower
replication rate $a_2$.  We suppose that when type 1
individuals replicate, the new individual has a type 1
genome with probability $1-r$ and a type 2 genome with
probability $r$. Type 2 genome individuals do not mutate.
The model is then
\begin{equation}\label{ode}
\begin{aligned}
{dv_1\over dt} &=a_1(1-r)v_1\\
{dv_2\over dt}&=a_1rv_1+a_2v_2 
\end{aligned}
\end{equation}
where $v_i$ is the number of type $i$ genomes for $i=1,2$. 
This is a variation of a model in Section 8.5 of
Nowak and May (2000).  A slightly different but perhaps
better interpretation of this model is to think of genome 1 as
being a specific (high performing) genome and genome 2 as
the collection of all the other genomes in the
population.     

This system
of differential equations is easily solved, and one can check that the ratio
$v_1/v_2$ converges as $t$ goes to infinity. It turns out that the limit is strictly positive
if and only if $r<1-a_2/a_1$. 
That is, in order for type 1 to be maintained in the population the mutation $r$ needs to
be below the threshold $1-a_2/a_1$. Hence, this model predicts that above a certain
mutation threshold faithful replication of the best adapted
genotype is compromised. Moreover, there seems to be general
agreement in the biology literature that above this
threshold the virus population will go extinct, see Eigen
(2002) and Manrubia et al.\ (2010). We propose here a simple
stochastic model that shows that this is not necessarily
so. In our model the population may survive, even if faithful
replication of the best adapted genotype is compromised, 
with the population being composed of ever changing
mutants. 

Our results may be biologically relevant for the following reason. 
An important current strategy to fight HIV
and other viruses is to try to increase the mutation
probability of the virus, see Eigen (2002) and Manrubia et
al (2010). This assumes that above a certain mutation
threshold the virus will die out. Our model suggests that at
least in theory this strategy may not work.

We now describe our continuous time evolution process.
Let $\mu$ be a probability
distribution with support contained in $[0,\infty)$ and
which is absolutely continuous with respect to Lebesgue
measure, and let $r\in[0,1]$.  Start with one individual at
time 0, and sample a birth rate $\la$ from the distribution
$\mu$. The individual gives birth at rate $\la$ and dies at
rate 1. Every time there is a birth the new individual:
(i) with probability $1-r$ keeps the same birth rate $\la$ as its parent,
and (ii) with probability $r$ is given a new birth rate $\la'$, sampled
independently of everything else from the
distribution $\mu$. 
We think of $r$ as the mutation probability and
the birth rate of an individual as representing the fitness
or genotype of the individual. Since $\mu$ is assumed
to be continuous, a genotype cannot appear more than once in
the evolution of the population. For convenience we label
the genotypes in the order of their appearance.

Let $Z(t)$ denote the number of individuals alive at time
$t$. We say that the evolution process \textit{survives} if
$Z(t)>0$ $\forall\ t\ge 0$ and otherwise \textit{dies out}. Our
main interest is in determining whether survival with 
positive probability is possible and by what
mechanism can survival be achieved.

\begin{thm} For $0\le r\le 1$ and probability distributions
  $\mu$ on $[0,\infty)$, the evolution process survives with
  positive probability if and only at least one of the
  following survival conditions holds:
  \[
\mu(\{\lambda: \lambda(1-r)>1\})>0,\leqno\emph{(I)}
\]
\[
\int_{\{\lambda:\lambda(1-r)\le 1\}}  \frac{\lambda r}{1-\lambda(1-r)}
d\mu(\lambda) >1 .
\leqno\emph{(II)}
\]
\end{thm}

The two extreme cases $r=0$ and $r=1$ are easy to
understand. If $r=0$ then (II) cannot hold and (I) reduces
to $\mu((1,\infty))>0$. In this case, no new types are ever
produced, the initial branching rate is used
forever by all individuals. Conditional on the initial branching rate
$\lambda$, $Z(t)$ is a linear birth-death process which
survives iff $\lambda>1$. Thus (I) is equivalent to positive
probability of survival. When
$r=1$, (I) cannot hold and (II) reduces to $\int \lambda
d\mu(\lambda)>1$. Now each new individual is a new genotype.
It is not hard to see that conditional on a given
individual's branching rate $\lambda$, the total number of
offspring of that individual is $k$ with probability 
\[
\frac{1}{1+\lambda}
\Bigl( \frac{\lambda}{1+\lambda}\Bigr)^k,
\quad k=0,1,\dots,
\]
with mean $\lambda$. 
Thus the unconditional
mean number of offspring of the first individual is
$\int \lambda\mu(d\lambda)$,
and the total number of individuals that ever live in
the evolutionary process is the same as the total progeny in
a Galton-Watson process with an offspring distribution which
has this mean. The total progeny is
infinite with positive probability if and only if this mean
is larger than 1, so (II) is equivalent to positive
probability of survival.

Condition (I) corresponds to the prediction of the
differential equation model \eqref{ode}. That is, below a
certain threshold for the mutation probability the virus can
survive because a well adapted (i.e. high $\la$) fixed
genotype can survive. However, if (I) fails it is still
possible to have survival by (II). In this case survival
holds because of a growing ``cloud'' of ever changing
mutants of low replicative ability.

Observe that for any $\epsilon>0$ and $r$ in $[0,1)$ 
there are distributions $\mu$ for which (I) holds but
$\int \lambda d\mu(\lambda)<\epsilon$.  This shows
that the behavior of our evolution process is drastically
different from the classical Galton-Watson process in
homogeneous or random environments. For these processes
survival is possible if and only if the expected
offspring (or a closely related expectation) is large
enough (see Harris (1989) for homogeneous
environments and Smith and Wilkinson (1969) for random
environments).

It is clear that if the support of $\mu$ is unbounded
then (I) holds for all $r<1$, so for interesting examples we
consider distributions with compact support. Among these 
distributions a natural family to consider is the 
uniform distribution on $[0,a]$,  $a>0$.
As the following shows, this class exhibits all possible
types of survival behavior
depending on the exact values of $a$ and $r$.

\begin{cor} Let $\mu$ be the uniform distribution on
  $[0,a]$, $a>0$. If $0<a\le 1$ then the evolution process
  dies out a.s.\ for all $r\in[0,1]$, while if $a>2$ the
  evolution process survives with positive probability for
  all $r\in[0,1]$.  If $a=2$ then the evolution process dies
  out a.s.\ for $r=1$ and survives with positive probability 
  for all $r\in[0,1)$. If $1<a<2$ then
  there exists $r_c\in(1-\frac1a,1)$
  such that \begin{itemize}
\item[\emph{(a)}] If $r<1-\frac1a$ then \emph{(I)} holds and the
  evolution process survives with positive probability.
\item[\emph{(b)}] If $1-\frac1a\le r<r_c$ then \emph{(II)} holds and 
the evolution process survives with positive probability.
\item[\emph{(c)}] If $r\ge r_c$ then the evolution process dies out a.s.
\end{itemize}
\end{cor}
\bigskip 

In words, whether the population goes extinct when the mutation rate is above a certain threshold depends crucially on the value of $a$. If $a>2$ there is no such threshold: the population survives for any mutation probability $r$. Note also that for $1<a<2$ there are two distinct thresholds: $1-\frac1a$ and $r_c$. If $r< 1-\frac1a$ a well adapted genome may survive forever while if $1-\frac1a\le r<r_c$ no fixed genome can survive forever. In this regime the population survives as a growing cloud of ever changing 
mutants. Finally, if $r\geq r_c$ the population goes extinct.

\section{Proof of Theorem 1}
\begin{proof}[Proof of Theorem 1]
  Recall that we start with a single genotype 1 individual
  at time 0.  Let $X_t$ be the number of type 1 individuals
  alive at time $t$.  Conditional on the initial branching
  rate $\la$, $X_t$ is a birth-death process with individual
  birth rate $\la(1-r)$ and death rate 1. In particular, it
  is well known (see Chapter 4 of Karlin and Taylor (1975))
  that it survives with positive probability if and only if
  $\la(1-r)>1$, and that
\begin{equation}\label{xtmean}
E(X_t|\la)=\exp((\la(1-r)-1)t).
\end{equation}
Integration of the condition $\lambda(1-r)>1$ with respect to
$\mu$ gives 
\[
P(X_t \geq 1\ \forall\ t>0)>0 
\text{ iff } 
\mu(\{\la:\la(1-r)>1\})>0.
\]

Now let $Y_t$ be the number of different genotypes born up
to time $t$ that are offspring of genotype 1 individuals.
Then $Y_t\uparrow Y_\infty$ as $t\to\infty$, the total
number of different genotypes ever produced by genotype 1
individuals. Note that if $r>0$ then $Y_\infty<\infty$ if and only if
$X_t=0$ eventually. For $h>0$ it is easy to see that
\[
E(Y_{t+h}-Y_t|\la,X_t)=\la r hX(t)+o(h) \text{ as
}h\downarrow 0,
\]
from which it follows that
$${d\over dt}E(Y_t|\la)=\la rE(X_t|\la)$$
and therefore, using \eqref{xtmean},
$$E(Y_t|\la)=r \la\int_0^t E(X_s|\la) ds=r\la\int_0^t \exp((\la(1-r)-1)s)ds .$$
Integration with respect to the measure $\mu$ now yields
\[
E(Y_t)=\int_0^{+\infty}\int_0^tr\la\exp((\la(1-r)-1)s)dsd\,\mu(\la).
\]
By the monotone convergence theorem, $E(Y_t)\uparrow
E(Y_\infty)$ as $t\to\infty$. Letting 
$\ds m(r)=E(Y_\infty)$, it is easy to show using the above
that 
\[
m(r)= \begin{cases}\ds
+\infty & \text{if }\mu(\{\lambda:\lambda(1-r)>1\})>0\\
\ds \\
\ds \int_0^{1/(1-r) }\frac{r \lambda}{1-\lambda(1-r)}d\mu(\lambda)
  & \text{if }\mu(\{\lambda:\lambda(1-r)>1\})=0 .
\end{cases}
\]

We now define the tree of genotypes first introduced by
Schinazi and Schweinsberg (2008) for a different model.
Assume (I) does not hold, and thus $Y_\infty <\infty$
a.s. Each vertex in the tree will be labeled by a positive
integer.  There will be a vertex labeled $k$ if and only if
an individual of genotype $k$ is born at some time.  We draw
a directed edge from $j$ to $k$ if the first individual of
genotype $k$ to be born had an individual of genotype $j$ as
its parent.  This construction gives a tree whose root is
labeled $1$ because all individuals are descendants of the
individual of genotype $1$ that is present at time zero.
The tree of genotypes is a (discrete time) Galton-Watson
tree with offspring distribution $p_k=P(Y_\infty=k)$. The
mean of the offspring distribution is $m(r)$, and hence, the
tree of genotypes is infinite with positive probability if
and only if $m(r)>1$.

To finish the proof, we claim that there are only two ways
for the evolution process to survive: either a fixed
genotype survives forever with positive probability
((I) holds), or the tree of genotypes is infinite
with positive probability ((II) holds). It is
clear that if either of these occur then the evolution
process survives with positive probability. Suppose now that
both (I) and (II) fail. Then with probability one each genotype that
ever appears gives birth to only finitely many individuals
and also the tree of types is finite a.s. This means that  
the total number of individuals that ever appear is
finite.
\end{proof}

\section{Proof of Corollary 1}
Let $\mu$ be the uniform distribution on $[0,a]$ where
$a>0$. Then (I) is equivalent to $a(1-r)>1$. 
If $a(1-r)\le 1$ then
\begin{equation}\label{mofr}
m(r) = \dfrac{1}{a} \int_0^a
\dfrac{r\lambda}{1-(1-r)\lambda} d\lambda .
\end{equation}

\paragraph*{The case $0<a\le 1$.}
Here $a(1-r)\le 1$ for all $r\in[0,1]$, so
(I) does not hold. Furthermore, the
fact that $a\le 1$ implies that the integrand in \eqref{mofr} is
an increasing function of $r$. Thus for all $r\in[0,1]$,
\[
m(r)\le m(1) =a/2 <1,
\]
and hence (II) also fails. For every $r$ the evolution
process dies out a.s.

\paragraph*{The case $a>1$.}
A little calculus shows that 
\begin{equation}\label{mofr2}
m(r)=-{r\over 1-r}-{1\over a}{r\over
  (1-r)^2}\ln(1-a(1-r)), \quad 
r \in (1-{1\over a},1).
\end{equation}
To complete the proof of Corollary 1 we will need the
following properties of $m(r)$. 

\begin{description}
\item[(P1)] $m(r)$ is continuous on $(1-1/a,1]$, 
$\ds \lim_{r\downarrow 1-1/a}m(r) = \infty\text{ and }
\lim_{r\uparrow 1}m(r) = a/2$.
\item[(P2)] If  $a\ge 3/2$ then $m(r)$ is strictly
  decreasing on  $(1-\frac1a,1)$ 
\item[(P3)]
If $1<a< 3/2$ then there exists $r_a\in
(1-\frac1a,1)$ such that $m(r)$ is strictly
decreasing on  $(1-\frac1a,r_a)$ and strictly increasing
  on $(r_a,1)$.
\end{description}
The proof of (P1) is simple and we will omit it. The proofs
of (P2) and (P3) require some work, so we will postpone them
for now and complete the proof of Corollary 1 assuming
(P2) and (P3) have been established. We consider three cases.

(i) If $a\ge 2$ and $r<1$, then by (P2) $m(r)> m(1)=a/2\ge1$,
so (II) holds for all $r\in(1-1/a,1)$. Also, $m(1)=a/2$
implies (II) holds for $a>2$ but fails for $a=2$.

(ii) If $3/2\le a<2$ then by $m(1)<1$, and hence by
(P1) and (P2) there
exists a unique $r_c\in(1-1/a,r_a)$ such that $m(r_c)=1$. By (P2),
(II) holds for $r<r_c$ but fails for $r\ge r_c$. 

(iii) If $1<a<3/2$ then by (P1) and (P3) $m(r_a)<1$. It follows
that there exists a unique $r_c\in(1-1/a,r_a)$ such that
$m(r_c)=1$, $m(r)>1$  
on $(1-1/a,r_c)$ and  $m(r)<1 $ on $(r_c,1]$. 

\bigskip The proof of Corollary 1 is now complete except
for the proofs of (P2) and (P3). 
At this point it is convenient to change variables. If we define
the function 
\[
g(x)=1-x+{1\over a}(x-x^2)\ln(1-{a\over x}), \quad x\in(a,\infty),
\]
then 
\[
m(r) = g(\frac1{1-r}) .
\]
Moreover, $m$ is increasing (decreasing) on the interval $(r_1,r_2)$ iff
$g$ is increasing (decreasing) on the interval $((1-r_1)^{-1},(1-r_2)^{-1})$.
A little calculation gives the first three derivatives of $g$,
\begin{align*}
g'(x) &=-1-{x-1\over x-a}-{1\over a}(2x-1)\ln(1-{a\over x})\\
g''(x)&=-{a-3ax+2x^2\over x(x-a)^2}-{2\over a}\ln(1-{a\over x})\\
g'''(x)&=-{a^2+ax(2a-3)\over x^2(x-a)^3}.
\end{align*}
With some additional calculation one can explicitly check
that
\begin{align}
\lim_{x\downarrow a}g'(x)=-\infty, \qquad
&\lim_{x\to +\infty}g'(x)=0,\label{g1}\\
\lim_{x\downarrow a}g''(x)=+\infty, \qquad
&\lim_{x\to +\infty}g''(x)=0.\label{g2}
\end{align}
We also note that by (P1),
\begin{equation}\label{g}
\lim_{x\downarrow a}g(x)=\infty, \qquad
\lim_{x\to +\infty}g(x)=a/2 .
\end{equation}

Suppose $a\ge
3/2$. Then  $g'''(x)<0$ for all $x>a$, and hence the
function $g''$ is strictly decreasing on
$(a,\infty)$. In view of \eqref{g2}, 
$g''$ must be positive on $(a,+\infty)$,
which implies $g'$ is strictly increasing on $(a,+\infty)$. 
In view of \eqref{g1}, $g'$ must be negative on $(a,+\infty)$, which
implies $g$ is strictly decreasing on  $(a,+\infty)$. This
means that $m(r)$ is strictly decreasing on $(1-1/a,1)$, so
(P2) is proved. 

Finally, suppose that $1<a<3/2$, and put $b=a/(3-2a)$.
Then $b>a$, $g'''<0$ on $(a,b)$ and $g'''>0$ on
$(b,\infty)$. As a consequence, $g''$ 
is strictly decreasing on $(a,b)$ and 
strictly increasing on $(b,\infty)$. In view of \eqref{g2}
there must exist a unique $c\in(a,b)$ such that
$g''>0$ on $(a,c)$ and $g''<0$ on $(c,\infty)$. This implies
$g'$ is strictly increasing on $(a,c)$ and strictly
decreasing on $(c,\infty)$. In view of \eqref{g1} 
there must exist a unique 
$x_a\in(a,c)$ such that $g'<0$ on $(a,x_a)$ and $g'>0$ on
$(x_a,\infty)$. This implies $g$ is strictly decreasing on
$(a,x_a)$ and strictly increasing on $(x_a,\infty)$. By
setting $r_a=1-1/x_a$ 
and using the correspondence between the functions $m$ and
$g$ we obtain (P3). 

\bigbreak
{\bf References.}
\medskip

M. Eigen (1971) Selforganization of matter and the evolution of biological macromolecules. Naturwissenschaften 58 465-523.

M. Eigen (2002) Error catastrophe and antiviral strategy. PNAS 99 13374-13376.

M. Eigen and P. Schuster (1977) The hypercycle. A principle of self-organization. 
Part A: emergence of the hypercycle.
Naturwissenschaften 64 541-565.

S. F. Elena and A. Moya (1999) Rate of deleterious mutation and the distribution of its effects on fitness in vesicular stomatitis virus. J. Evol. Biol. 12 1078-1088.

T.E.Harris (1989) {\it The Theory of Branching Processes.} Dover, New York. 

S. Karlin and H.M. Taylor (1975) {\it A First Course in
  Stochastic Processes, 2nd Ed.} Academic Press, New York.

S.C. Manrubia, E. Domingo and E. Lazaro (2010) Pathways to extinction: beyond the error threshold. Phil. Trans. R. Soc. B 365, 1943-1952.

M. A. Nowak and R.M.May (2000) {\it Virus dynamics.}  Oxford University Press.

R. Sanjuan, A. Moya and S. F. Elena (2004). The distribution of fitness effects caused by single-nucleotide substitutions in an RNA virus. PNAS 1018396-8401.

R.B.Schinazi and J. Schweinsberg (2008) Spatial and non spatial stochastic models for immune response. Markov Processes and Related Fields 14 255-276.

W.L.Smith and W.E.Wilkinson (1969) On branching processes in random environments. The Annals of Mathematical Statistics 40, 814-827.

M. Vignuzzi, J.K. Stone, J.J. Arnold, C.E.Cameron and R. Andino (2006). Quasispecies diversity determines pathogenesis through cooperative interactions in a viral population. Nature 439, 344-348

\end{document}